\documentclass[11pt]{article}
\usepackage[utf8]{inputenc}

\usepackage[bottom=4cm, left =2.5cm, right =2.5cm]{geometry} 
\usepackage{import} 
\usepackage{fancyhdr}
\usepackage{hyperref} 
\usepackage{amsthm, thmtools,mathtools} 
\usepackage[intoc, english]{nomencl} 

\usepackage{amsmath, amsfonts, amssymb} 
\usepackage{stmaryrd} 
\usepackage{colonequals} 
\usepackage{bbm} 
\usepackage{enumitem} 
\usepackage[normalem]{ulem} 

\usepackage{lmodern} 

\usepackage{pgf} 
\usepackage{graphicx} 
\usepackage{float} 
\usepackage{xcolor} 

\allowdisplaybreaks[1]
\usepackage[english]{babel} 
\usepackage[T1]{fontenc} 
\usepackage{authblk} 

\usepackage{comment}

\usepackage[textwidth=2.0cm, textsize=tiny]{todonotes}

\usepackage[backend=biber,style=alphabetic,sorting=nyt,giveninits=true,maxnames=99]{biblatex}
\usepackage{csquotes}
\newtheorem{thm}{Theorem}[section]
\newtheorem{prp}[thm]{Proposition}
\newtheorem{lem}[thm]{Lemma}
\newtheorem{cor}[thm]{Corollary}

\newtheorem{df}[thm]{Definition}
\newtheorem{prp-df}[thm]{Definition-Proposition}

\theoremstyle{definition}
\newtheorem{nota}[thm]{Notations}

\theoremstyle{remark}
\newtheorem{rem}[thm]{Remark}

\newtheorem{exm}[thm]{Example}

\newcommand{\comments}[1]{~}
\newcommand{\E}{\mathbb{E}}
\newcommand{\condEp}[2]{\left.\E\left(#1\right\lvert #2 \right)}
\newcommand{\condEc}[2]{\E\left.\left[#1\right\lvert #2 \right]}
\newcommand{\R}{\mathbb{R}}
\newcommand{\Proba}{\mathbb{P}}
\newcommand{\law}{\mathcal{L}}

\newcommand{\N}{\mathbb{N}}
\newcommand{\Z}{\mathbb{Z}}

\newcommand{\calF}{\mathcal{F}}

\newcommand{\indic}{\mathbbm{1}}
\newcommand{\norm}[1]{\left\Vert #1\right\Vert}
\newcommand{\normB}[1]{\norm{#1}_{\mathbb B}}
\newcommand{\abs}[1]{\left\lvert #1 \right\rvert}

\newcommand{\eqdef}{\colonequals}

\newcommand\bigp[1]{\left(#1\right)}

\newcommand{\bigcro}[1]{\left[#1\right]}

\newcommand{\eps}{\varepsilon}

\renewcommand{\leq}{\leqslant}
\renewcommand{\geq}{\geqslant}

\renewcommand{\Tilde}{\widetilde}

\newcommand{\finpreuve}{\\{\color{white}.}\hspace{\linewidth}\hspace{-12pt} $\qedsymbol$\\}

\definecolor{DarkBlue}{rgb}{0.0, 0.28, 0.39}

\makeatletter
\@addtoreset{equation}{section}
\makeatother

\hypersetup{
    colorlinks=true,       
    linkbordercolor=DarkBlue,
    linkcolor=DarkBlue!70!DarkBlue, 
    citecolor=DarkBlue!70!DarkBlue, 
    urlcolor=DarkBlue!70!DarkBlue 
}

\AtBeginBibliography{\small}

\begin{document}
    \title{Central limit theorem under the Dedecker-Rio condition in some Banach spaces}
    \author{Aurélie Bigot\thanks{LAMA, Univ Gustave Eiffel, Univ Paris Est Créteil, UMR 8050 CNRS, F-77454 Marne-La-Vallée, France.}}
    \date{}
    \maketitle

\begin{abstract} 
    We extend the central limit theorem under the Dedecker-Rio condition to adapted stationary and ergodic sequences of random variables taking values in a class of smooth Banach spaces. This result applies to the case of random variables taking values in $L^p(\mu)$, with $2 \leq p < \infty$ and $\mu$ a $\sigma$-finite real measure. As an application we give a sufficient condition for empirical processes indexed by Sobolev balls to satisfy the central limit theorem, and discuss about the optimality of these conditions. 
\end{abstract}

\section*{Introduction}
~\vspace{-15pt}
\\

Let $(\Omega, \mathcal F, \Proba)$ be a probability space and $(X_i)_{i \in \mathbb Z}$ be a strictly stationary sequence of centered and square integrable real-valued random variables, adapted to a stationary filtration $(\calF _i)_{i \in \mathbb Z}$ . In 2000, Dedecker and Rio proved in \cite{DR2000} the central limit theorem (in short CLT) for $(X_n)_n$ under the condition 
\begin{align*}\label{conditionDR}
    \tag{\textbf{\textit{DR}}}
    X_0 \condEp{S_n}{\calF_0} \quad \textrm{converges in } \mathbb L ^1
\end{align*}
where $S_n = X_1 + \cdots + X_n$. 
Recall that for $\alpha$-mixing random variables, condition (\ref{conditionDR}) leads to the following condition:
\begin{align}
        \label{conditionDR_alpha}
    \sum_{n =0}^{+ \infty} \int_0^{\alpha(\calF_0, \sigma(X_n))} Q_{\abs{X_0}}^2 (u) du \, < \infty,
\end{align}
where $Q_{\abs{X_0}}$ is the upper tail quantile function of $ \abs{X_0}$, and $(\alpha(\calF _0, \sigma(X_n)))_{n \geq 1}$ is the sequence of strong mixing coefficients associated with $(X_i)_{i \in \mathbb Z}$ (see for instance \parencite[(2.1)]{DR2000}). Note that condition (\ref{conditionDR_alpha}) is known to be essentially optimal as proved in \parencite[Section 4]{DMR1994} and in \cite{Bra1997}. 
\\

The CLT under condition (\ref{conditionDR}) has been extended in \cite{DM2003} to random variables taking values in a separable Hilbert space. In this paper we extend this CLT to the case of r.v.'s taking values in a 2-smooth Banach space with a Schauder basis. As we shall see, the main ingredients to obtain such a result are a martingale blocks decomposition,  Theorem 5 in \parencite[]{Ros1982} and Theorem 2.1 in \cite{DM2015}. Typically, the $L^p$ spaces for $p \geq 2$ fit into our framework of Banach spaces, and this particular case will lead to a CLT for the empirical process as will be shown in Section 2. 

Recently, several authors have extended other projective criteria valid in the real spaces context to the case of smooth Banach spaces.
For instance, the Hannan condition (see \cite{Han1973}) has been extended in \cite{DMP2013} for random variables taking values in a 2-smooth Banach space having a Schauder basis. Such a condition can be written in the Banach space setting: $\sum_{k \in \Z} (\E\norm{P_0(X_k)}_{\mathbb B}^2)^{1/2} < \infty$ where $P_0$ is the operator defined by $P_0 = \condEp{\cdot}{\calF_0}-\condEp{\cdot}{\calF_{-1}}$, $\mathbb B$ is the real and separable Banach space and $\normB{\cdot}$ its associated norm. 
Very recently the condition of Maxwell-Woodroofe (see \cite{MW2000}) has been extended by \cite{Cun2017} in Banach space settings. In the case of 2-smooth Banach spaces, the condition becomes $\sum_{n = 1}^\infty n^{-3/2}(\E[\norm{\condEp{S_n}{\calF_0}}_{\mathbb B}^2])^{1/2} < \infty$.  
Note that it has been shown in \cite{DV2008} that the conditions by Dedecker-Rio, Hannan and Maxwell-Woodroofe are independent. 

The paper is organised as follows. In Section 1, we state an extension of the CLT under the condition (\ref{conditionDR}) to the case of 2-smooth Banach spaces with a Schauder basis. As a consequence, in Section 2, we derive in Corollary \ref{Lp_thm} a sufficient condition for empirical processes indexed by Sobolev balls to satisfy the CLT, and discuss about the optimality of this condition. The proofs of the main results are postponed to Section 3.  

\section{A CLT in some smooth Banach spaces}
~
In all the paper, $(\mathbb B, \norm{\cdot}_{\mathbb B})$ will be a real and separable Banach space. We shall consider the class of Banach spaces that are 2-smooth. This notion introduced by Pisier in \parencite[Section 3]{Pis1975} plays the same role with respect to vector martingales as spaces of type 2 do with respect to the sums of independent random vectors. Let us consider the following definition of 2-smooth Banach spaces.

\begin{df}
    Let $(\mathbb B, \norm{\cdot}_{\mathbb B})$ be a separable Banach space and define $\psi _2 : x \mapsto \norm{x}_{\mathbb B}^2$. $(\mathbb B, \norm{\cdot}_{\mathbb B})$ is said to be 2-smooth if there exists $d > 0$ such that for any $x,u \in \mathbb B$
    \begin{enumerate}[label = (\roman*)]
        \item  if $x \neq 0$, $D^2\psi_2(x)(u,u) \leq d^2 \norm{u} _ {\mathbb B}^2$
        \item $\abs{D\psi_2(x)(u) - D\psi_2(0)(u)} \leq d^2\norm{x}_{\mathbb B}\norm{u}_{\mathbb B}$.
    \end{enumerate}
    
    Here $D\psi_2(x)$ and $D^2\psi_2(x)$ denote respectively the usual first and second order Fréchet derivative of $\psi_2$ at point $x$.
\end{df}

As quoted in \cite{Pin1994}, this definition implies the $(2,d)$-smoothness in the sense of Pisier meaning that 
\begin{align*}
    \norm{x+y}_{\mathbb B}^2 + \norm{x - y}_{\mathbb B}^2 \leq 2 \norm{x}_{\mathbb B}^2 + 2d^2\norm{y}_{\mathbb B}^2, \;\; \forall x,y \in \mathbb B
    .
\end{align*}

We shall also need the notion of Banach spaces with a Schauder basis. 
\begin{df}
    A family $\{x_n \, : \, n \in \N \}$ of elements of $\mathbb B$ is called a Schauder basis if for any vector $x \in \mathbb B$ there exists an unique series 
    $
        \sum_{n=0}^\infty a_n x_n $, with $a_n = a_n(x) \in \R
    $, 
    which converges to $x$ with respect to $\norm{\cdot}_{\mathbb B}$. 
\end{df}

One of the properties of Banach spaces that admit a Schauder basis is the uniform boundedness of the family of operators $(P_n)_n$, where $P_n$ is the projection on the space generated by $\{x_k \, : \, k \leq n\}$. More precisely, there exists a constant $c>0$ such that for any $n \in \N$, $\norm{P_n} \leq c$, where $\norm{\cdot}$ is the operator norm.
\\

\begin{exm}
    \label{exm:Lp}
    Let $p \in [1, +\infty[$ and $\mu$ a $\sigma$-finite measure on $\R$. 
    According to \parencite[Section 6]{AK2016}, $L^p(\mu)$ is isometrically isomorphic to one of the following spaces: 
    \begin{align*}
    (\mathbb R^n, \norm{\cdot}_{\ell^p}), \ell^p,  L^p(\nu), L^p(\nu) \oplus (\mathbb R^n, \norm{\cdot}_{\ell^p}), L^p(\nu) \oplus \ell^p
    \end{align*}
    where $\nu$ is a nonatomic probability measure on a countably generated measurable space. 
    On the one hand, in $(\R^n, \norm{\cdot}_{\ell^p})$ the canonical basis is a Schauder basis and in $\ell^p$ the canonical basis is a Schauder basis. On the other hand, $L^p(\nu)$ is isometrically isomorphic to $L^p([0,1], \lambda)$ where $\lambda$ is the Lebesgue measure, and according to \parencite[Proposition 6.1.3]{AK2016}, $L^p([0,1], \lambda)$ has a Schauder basis. Hence, for any $p \in [1, +\infty[$, $L^p(\mu)$ is a Banach space with a Schauder basis. 
    \\
    In addition, equipped with its usual norm, for any $p\geq 2$, $L^p(\mu)$ is $(2, \sqrt{p-1})$-smooth as shown in \parencite[Proposition 2.1]{Pin1994}. 
\end{exm}

The main result of this paper is Theorem \ref{thm_ppal} below which extends \parencite[Theorem 1]{DR2000} to the Banach space setting. 

\begin{thm}\label{thm_ppal}
    Let $(\mathbb B, \norm{\cdot}_{\mathbb B})$ be a 2-smooth Banach space with a Schauder basis. Let $(X_i)_{i \in \Z}$ be an ergodic stationary sequence of centered $\mathbb B$-valued random variables, adapted to a non-decreasing and stationary filtration $(\mathcal F _i)_{i \in \Z}$ and such that $\E (\norm{X_0}_{\mathbb B}^2) < \infty$. 
    Set $S_n \eqdef \sum\limits_{k = 1}^n X_k$ and assume that 
    \begin{align}\label{thm_ppal:H}
        \norm{X_0}_{\mathbb B} \condEp{S_n}{\calF_0}  \text{ converges in } \mathbb{L}^1_{\mathbb B}.
    \end{align}
    Then $\bigp{\frac{1}{\sqrt n}S_n}_{n \geq 1}$ converges in distribution to $G$, where $G$  is a Gaussian $\mathbb B$-valued random variable whose covariance operator is given by: $K_G(x^*, y^*) = \sum_{k \in \Z} cov(x^*(X_0), y^*(X_k))$ for any $x^*, y^* \in \mathbb B ^*$. 
\end{thm}

In view of applications, we shall give in Corollary \ref{cor_cond_suff} sufficient conditions for (\ref{thm_ppal:H}) to be verified. With this aim, we first introduce useful notations.

\begin{nota}
    Let $X$ and $Y$ be real-valued random variables. Denote by : 
    \begin{itemize}
        \item $Q_{\abs{X}}$ the generalized inverse of the upper tail function $t \mapsto \Proba(\abs{X} >t)$
        \item $G_{\abs{X}}$ the inverse of $x \in [0, \Proba(\abs X >0)] \mapsto \int_0^x Q_{\abs{X}}(u) du$
        \item $H_{X,Y}$ the generalized inverse of $x \mapsto \E(\abs{X} \indic_{\abs{Y}>x})$. 
    \end{itemize}    
\end{nota}

\comments{Following the notion introduced by Volkonskii and Rozanov in \cite{RV1959}, we define the classical $\beta$-mixing coefficients.}
\begin{df}
    Consider a stationary sequence of random variables $(X_i)_{i \in \mathbb Z}$ adapted to a non-decreasing and stationary filtration $(\calF _i)_{i \in \Z}$. We define for any nonnegative integer $k$ : 
    \begin{align}
        \beta_{1,X}(k) = \norm{ \sup_{\norm{f}_\infty \leq 1} \abs{ P_{X_k \lvert \mathcal F _0}(f) - P_{X_k}(f) } }_1
        .
        \nonumber
    \end{align}
\end{df}
If there is no confusion on the r.v. to which we refer, we shall write $\beta_1(k)$ instead of $\beta_{1,X}(k)$ for the sake of clarity. 

\begin{cor}\label{cor_cond_suff}
    Let $(\mathbb B, \norm{\cdot}_{\mathbb B})$ be a 2-smooth Banach space with a Schauder basis. Let $(X_i)_{i \in \Z}$ be an ergodic stationary sequence of centered $\mathbb B$-valued random variables, adapted to a non-decreasing and stationary filtration $(\mathcal F _i)_{i \in \Z}$ and such that $\E (\norm{X_0}_{\mathbb B}^2) < \infty$. 
    Set $\gamma_i \eqdef \E(\norm{\condEc{X_i}{\calF_0}}_{\mathbb B})$. Consider the conditions 
    \begin{enumerate}[label = (\roman*)]

        \item $\sum_{n \geq 1} \int_0^{\beta_1(n)} Q_{\norm{X_0}_{\mathbb B}}^2(u) \, du \, < \infty$,
        \label{cor_cond_suff:item1}
        
        \item $\sum_{n \geq 1} \int_0^{\gamma_n} Q_{\norm{X_0}_{\mathbb B}} \circ G_{\norm{X_0}_{\mathbb B}} (u) \, du \, < \infty$ 
        \label{cor_cond_suff:item2}.
    \end{enumerate}
    We have the implications $(i) \Rightarrow (ii) 
    \Rightarrow$ (\ref{thm_ppal:H}). 
\end{cor}
~\\
\textit{Proof of Corollary \ref{cor_cond_suff}.}
    By Berbee's coupling lemma (see \parencite[Lemma 5.1]{Rio2017}), there exists a random variable $X_n^*$ distributed as $X_n$, independent of $\calF_0$ and such that $\Proba(X_n \neq X_n^*) = \beta_1(n)$.
    Hence
    \begin{align}
        \gamma_n 
        = \E\bigp{ \norm{\condEp{X_n - X_n^*}{\calF_0} }_{\mathbb B }}
         &\leq \E \bigp{ \norm{X_n-X_n^*}_{\mathbb B} \indic_{X_n \neq X_n^*}}
        \nonumber
        \\
        & \leq \E \bigp{ \norm{X_n}_{\mathbb B} \indic_{X_n \neq X_n^*}} + \E \bigp{ \norm{X_n^*}_{\mathbb B} \indic_{X_n \neq X_n^*}}.
        \nonumber
    \end{align}
    Using the same argument as in the proof of Proposition 1 in \cite{DD2003}, we get 
    \begin{align*}
        \E \bigp{ \norm{X_n}_{\mathbb B} \indic_{X_n \neq X_n^*}} 
        \leq 
        \int_0^{\Proba(X_n \neq X_n^*)} Q_{\norm{X_n}_{\mathbb B}}(u) du,
    \end{align*}
    implying that 
    \begin{align}
        \gamma_n 
        & \leq 2 \int_0^{\beta_1 (n)} Q_{\norm{X_n}_{\mathbb B}}(u) du 
        \quad \textrm{and then} \quad G_{\norm{X_0}_{\mathbb B}}^{-1}(\gamma_n/2) \leq \beta_1(n).
        \nonumber
    \end{align}
    Therefore, using a change of variables, it follows that
    \begin{align}
        \int_0^{\gamma_n} Q_{\norm{X_0}_{\mathbb B}} \circ G_{\norm{X_0}_{\mathbb B}} (u) \, du
        \leq 
        2\int_0^{\gamma_n /2} Q_{\norm{X_0}_{\mathbb B}} \circ G_{\norm{X_0}_{\mathbb B}} (u) \, du
        \leq 
        2\int_0^{\beta_1(n)} Q^2_{\norm{X_0}_{\mathbb B}}(u) \, du
        .
        \nonumber
    \end{align}
    ~
    \\
    This ends the proof of the first implication. The second implication is an immediate consequence of the proof of Proposition 1 in \cite{DD2003} by replacing the absolute values by the norm $\norm{\cdot}_{\mathbb B}$. 
    \comments{We turn to the proof of the second implication. 
    \\
    Following the proof of Proposition 1 in \cite{DD2003}, for any $n$,
    \begin{align}
        \E\bigcro{ \norm{\condEp{X_n}{\calF_0}}_{\mathbb B} \indic_{\norm{X_0}_{\mathbb B}> t } } 
        \leq 
        \E( \norm{\condEp{X_n}{\calF_0}}_{\mathbb B} ) \wedge \E\bigp{ \norm{X_n}_{\mathbb B} \indic_{\norm{X_0}_{\mathbb B} > t } }
        \nonumber
        \\
        \leq 
        \int_0^{\gamma_n} \indic_{u < \E\bigp{ \norm{X_n}_{\mathbb B} \indic_{\norm{X_0}_{\mathbb B} > t } } } du 
    \end{align}
    so that using Tonelli's theorem :
    \begin{align}
        \E\bigcro{ \norm{\condEp{X_n}{\calF_0}}_{\mathbb B} \norm{X_0}_{\mathbb B} }
        & = 
        \int_0^\infty \E\bigcro{ \norm{\condEp{X_n}{\calF_0}}_{\mathbb B} \indic_{\norm{X_0}_{\mathbb B}> t } }\, dt
        \nonumber
        \\
        & \quad \leq 
        \int_0^{\gamma_n} H_{\norm{X_n}_{\mathbb B}, \norm{X_0}_{\mathbb B}}(u) \, du
        \nonumber
        \\
        & \qquad \quad \leq 
        \int_0^{\gamma_n} Q_{\norm{X_0}_{\mathbb B}} \circ G_{\norm{X_0}_{\mathbb B}} (u) \, du
        .
        \nonumber
    \end{align}
    The last inequality follows from the Fréchet's inequality(1957) which implies in particular that : \\$\E(\abs X \indic_{\abs Y > t}) \leq \int_0^{\Proba (\abs Y > t)} Q_X(u) du$. } 
    \finpreuve

In view of applications, let us give the following result which specifies the rates of decrease of $(\beta_1(k))_{k > 0}$ and moments of $\norm{X_0}_{\mathbb B}$ for (\ref{thm_ppal:H}) to hold. Its proof follows directly from \parencite[Annex C with $p = 2$]{Rio2017}. 

\begin{cor}\label{cor_suff_cond}
    Let $(\mathbb B, \norm{\cdot}_{\mathbb B})$ be a 2-smooth Banach space with a Schauder basis. Let $(X_i)_{i \in \mathbb Z}$ be an ergodic stationary sequence of centered $\mathbb B$-valued random variables such that $\E \norm{X_0}_{\mathbb B}^2 < \infty$, and adapted to a non-decreasing and stationary filtration $(\mathcal F _k)_k$. Assume that one of the following conditions holds:
    \begin{enumerate}[label = (\roman*)]
        \item there exists $r > 2$ such that $\E(\norm{X_0}_{\mathbb B}^{r}) < \infty$ and $\sum_{n \geq 0} (n+1)^{2/(r-2)}  \beta_{1}(n) < \infty$

        \item there exist $r > 2$ and $c > 0$ such that for any $x$, $\Proba\bigp{ \norm{X_0}_{\mathbb B}^{} > x} \leq \bigp{\frac{c}{x}}^r$ and $\sum_{n \geq 0} \beta_{1}(n)^{1-2/r} < \infty$

        \item there exist $a > 0$ and $\tau > 0$ such that $\E \bigcro{\norm{X_0}_{\mathbb B}^{2} (\log (1+\norm{X_0}_{\mathbb B}^{}))^a } < \infty$ and $\beta_{1}(n) = O \bigp{e^{-\tau n^{1/a}}}$.
    \end{enumerate}
    Then $\sum_{n \geq 1} \int_0^{\beta_1(n)} Q_{\norm{X_0}_{\mathbb B}}^2(u) \, du \, < \infty$  is verified and the conclusion of Theorem \ref{thm_ppal} holds. 
\end{cor}

\section{Applications to the empirical processes over Sobolev balls for dependent sequences}

~
    Let us consider $(Y_i)_{i \in \Z}$ a stationary and ergodic sequence of real random variables, whose cumulative distribution function is denoted $F$, and define $F_n(t) = \frac 1 n \sum_{k = 1}^n \indic_{Y_k \leq t}$ the empirical distribution function. We are interested in the asymptotic behavior of the centered empirical distribution function in $L^p(\mu)$, where $p\geq 2$ and $\mu$ is a $\sigma$-finite measure on $\mathbb R$.
\\
    We suppose that
    \begin{align}
        \int_{\R_-} F(t)^p \, d\mu(t) + \int_{\R_+} (1-F(t))^p \, d\mu(t) < \infty
        \label{Lp_condition_F}
    \end{align}
    so that $F_n - F$ is a random element of $L^p(\mu)$. 
\\
    In \cite{DM2007} the link between the convergence in distribution of $\sqrt n (F_n - F)$ in $ L ^p(\mu)$ and Donsker classes has been clearly established. More precisely, let us denote 
    \begin{align*}
        W_{1,q}(\mu) \eqdef \left\{f : \R \rightarrow \R \, : \, f(x) = f(0) + \indic_{x >0} \int_{[0,x[}g \, d\mu \, - \indic_{x\leq 0} \int_{[x, 0[} g \, d\mu \, , \, \norm{g}_{q, \mu} \leq 1 \right\}
        ,
    \end{align*}
    where $q$ is the conjugate exponent of $p$ and $\norm{\cdot}_{q,\mu}$ is the usual norm on $L^q(\mu)$. Then according to \parencite[Lemma 1]{DM2007} the following convergences are equivalent: 
    \begin{enumerate}[label = (\roman*)]
            \item $
            \{\sqrt n (F_n - F)(t) \}_{_t} \xrightarrow[]{\law} \{ G(t)\}_{_t} \, \textrm{in } L^p(\mu)
            $
    
            \item $
            \left\{ \sqrt n \bigp{\frac 1 n \sum_{k= 1}^n f(Y_k) - \E f(Y_0)} \right\}
            \xrightarrow[]{\law} \{G_1(f)\} \, \textrm{in } \ell^\infty (W_{1,q}(\mu))
            $   
    \end{enumerate}
    where $\ell^\infty (W_{1,q}(\mu))$ is the space of all functions $\phi : W_{1,q}(\mu) \to \R$ such that $\sup_{f \in W_{1,q}(\mu)} \abs{\phi(f)}$ is finite, and $G_1(f) = \int g(t)G(t)d\mu(t)$. 
    \\
    Hence, proving that $W_{1,q}(\mu)$ is a Donsker class for $(Y_i)_{i \in \Z}$ is equivalent to proving that $\sqrt n (F_n - F)$ converges weakly in $L^p(\mu)$ to a Gaussian process. 
    \\
    
    \noindent 
    To study the asymptotic behavior of $\sqrt n (F_n - F)$, we then define the random process: 
    \begin{align}
        \forall i \in \Z, \, X_i = \{ \indic_{Y_i \leq t} - F(t)\, : \, t \in \R\}
        \label{Lp_def_X}
    \end{align}
    which takes values in $L^p(\mu)$. With such a notation, the study of the asymptotic behavior of $\bigp{\sqrt n (F_n - F)}_{n\geq 1}$ is equivalent to the study of the asymptotic behavior of $(S_n/\sqrt n)_{n \geq 1}$ in $L^p(\mu)$ where $S_n = X_1 + \cdots + X_n$. Hence, since $L^p(\mu)$, $p\geq 2$, is a 2-smooth Banach space with a Schauder basis, the centered empirical distribution behavior in $L^p(\mu)$ will follow from an application of Theorem \ref{thm_ppal} and in particular of Corollary \ref{cor_cond_suff}. 
    To state the condition in terms of dependence conditions on the sequence $(Y_i)_{i \in \Z}$, we introduce the following weak dependence coefficients (see \cite{DP2005}).
\begin{df}\label{def_weak_beta}
    Let $(Y_i)_{i \in \mathbb Z}$ be a stationary sequence of real random variables adapted to a stationary filtration $(\calF_i)_{i \in \mathbb Z}$. 
    For any nonnegative integer $k$, let  
    \begin{align*}
        b_k
        = \sup\limits_{t \in \R} \abs{ \condEp{\indic_{Y_k \leq t}}{\mathcal M_0} - \E(\indic_{Y_k \leq t}) }
        \textrm{ and }
        \Tilde \beta_{1,Y}(k) \eqdef 
        \E(b_k)
        .
    \end{align*}
\end{df}
\noindent
Let us consider the same notations as in \cite{DM2007}: 
\begin{nota}
    Define the function $F_\mu$ by: $F_\mu(x) = \mu(]0,x])$ if $x \geq 0$ and $F_\mu(x) = -\mu([x,0[)$ if $x \leq 0$. Define
also the nonnegative random variable $Y_{p,\mu} = \abs{F_\mu (Y_0)}^{1/p}$. 
\end{nota}
\noindent
As an application of Corollary \ref{cor_cond_suff}, we derive the following result. 
\begin{cor}\label{Lp_thm}
    Let $p \in [2, +\infty[$ and assume that 
    \begin{align}\label{Lp_thm:H}
        \sum_{n \geq 0} \int_0^{\Tilde \beta_{1,Y}(n)} Q_{Y_{p, \mu}}^2 (u) \, du \, < \infty.
    \end{align}
    Then 
    \begin{align}\label{Lp_thm:conclusion}
        \sqrt n (F_n - F) \xrightarrow[n \to +\infty]{\law} G \quad \textrm{in } L^p(\mu),
    \end{align}
    where $G$ is a Gaussian process whose covariance operator is given, for any $x^*, y^* \in L^q(\mu)$, by 
    $K_G(x^*, y^*) = \sum_{k \in \mathbb Z} cov(x^* (X_0), x^*(X_k))$. 
    In particular, 
    \begin{align}
        \label{Lp_thm:conclusion_reformul}
        n^{p/2} \int_{\R} \abs{F_n(t) - F(t)}^p \, d\mu(t)  \xrightarrow[n \to +\infty]{\law} \int_{\R} \abs{G(t)}^p \, d\mu(t)
        .
    \end{align}
\end{cor}

\begin{rem}\label{Lp_thm:H_case_finite}
    When $\mu$ is finite, (\ref{Lp_thm:H}) simply reads as $\sum_{n \geq 0} \Tilde \beta_{1,Y}(n) < \infty$. 
\end{rem}

With the equivalence with Donsker classes in mind, it is interesting to study the envelope function of our class of functions.
Note that $x \mapsto \abs{F_\mu(x)}^{1/p}$ is the smallest envelope function of $W_{1,q}(\mu)$ which means that
\begin{align}\label{def_smallest_env_fx}
    \abs{F_\mu(x)}^{1/p} = \sup\{\abs{f(x) - f(0)} \, : \, f \in W_{1,q}(\mu)\}
    .
\end{align}
Indeed, for any $f \in W_{1,q}(\mu)$ and any real $x$, by Hölder's inequality,
\begin{align}
    \abs{f(x) - f(0)} \leq \abs{F_\mu (x)}^{1/p} \norm{g}_{q, \mu} \leq \abs{F_\mu (x)}^{1/p}
    .
    \nonumber
\end{align}
On the other hand, define for any $x \in \R$: 
\begin{align}
    &f_x : y \mapsto \indic_{y >0} \int_{[0,y[}g_x \, d\mu \, - \indic_{y\leq 0} \int_{[y, 0[} g_x \, d\mu 
    \nonumber
    \\
    &\hspace{5cm} \textrm{ with } 
    g_x : t \mapsto
    \left\{
    \begin{array}{ll}
       \indic_{[0,x]}(t)\, \abs{F_\mu(x)}^{-1/q}  & \textrm{if } F_\mu(x) \neq 0 \textrm{ and } x>0 \\
       -\indic_{[x,0]}(t)\, \abs{F_\mu(x)}^{-1/q}  & \textrm{if } F_\mu(x) \neq 0 \textrm{ and } x<0 \\
       0  & \textrm{either}
    \end{array}
     \right. .
    \nonumber
\end{align}
One has that $f_x$ belongs to $W_{1,q}(\mu)$ for any real $x$ and $f_x(x) = \abs{F_\mu(x)}^{1/p}$. This ends the proof of (\ref{def_smallest_env_fx}). 
\\

\noindent
Note that when $(Y_i)_{i \in \mathbb Z}$ is a sequence of i.i.d. random variables, (\ref{Lp_thm:H}) reads as
\begin{align}
    \int_0^1 Q_{Y_{p,\mu}}^2(u) \, du < \infty \qquad \textrm{i.e. } \E[\abs{F_\mu(Y_0)}^{2/p}] < \infty.
    \nonumber
\end{align}
It means that the smallest envelope function of $W_{1,q}(\mu)$ is square integrable. 
This condition together with an entropy condition ensure the CLT (see Theorem 2.5.2 in \cite{VW1996}). 
\\ 

\begin{rem}
    It is worth noting that, in the particular case where $\mu$ is the Lebesgue measure on $\R$ denoted by $\lambda$, $W_{1,q}(\lambda)$ can be rewritten
    \begin{align*}
         W_{1,q}(\lambda) = \left\{f : \R \rightarrow \R \, : \, f(x) = f(0) + \indic_{x >0} \int_{[0,x[} g(t) \, dt \, - \indic_{x\leq 0} \int_{[x, 0[} g(t) \, dt \, , \, \norm{g}_{q} \leq 1 \right\} 
    \end{align*}
    which is the space of absolutely continuous functions $f$ such that $\lambda(\abs{f'}^q) \leq 1$. In particular, it contains the unit ball of the Sobolev space of order 1 with respect to $L^q(\lambda)$. 
    Moreover, since in this case $F_\mu$ is the identity function, condition 
    (\ref{Lp_thm:H}) can be rewritten 
    \begin{align*}
        \sum_{n \geq 0} \int_0^{\Tilde \beta_{1,Y}(n)} Q_{\abs{Y_0}}^{2/p} (u) \, du \, < \infty.
    \end{align*}
\end{rem}
~\\

\noindent
\textit{Comment on the optimality of condition (\ref{Lp_thm:H}).} \\
In the dependent case, condition (\ref{Lp_thm:H}) of Corollary \ref{Lp_thm} implies that 
\begin{align}\label{Lp_envelope_fx_convergence}
    \frac{1}{\sqrt n}\sum_{k=1}^n (|F_\mu(Y_k)|^{1/p}-\E(|F_\mu(Y_k)|^{1/p}))
\end{align}
converges in distribution to a Gaussian r.v.
The proposition below which is essentially due to \parencite[]{DMR1994}, shows that condition (\ref{Lp_thm:H}) is essentially optimal for the convergence in distribution of (\ref{Lp_envelope_fx_convergence}) to a Gaussian r.v. to hold. 

\begin{prp}[Doukhan, Massart, Rio]\label{Lp_thm:optimality}
    Let $a > 1$. Suppose that $Y_0$ is a real-valued r.v. whose distribution function $F$ is continuous and such that $\int_0^1 u^{-1/a} Q_{Y_{p,\mu}}^2(u) \, du = +\infty$. Then, there exists a stationary Markov chain $(Z_i)_{i \in \Z}$ with marginal distribution function $F$ and such that
    \begin{enumerate}[label = (\roman*)]
        \item $0 < \liminf\limits_{n \to \infty} n^a \beta_1(n) \leq \limsup\limits_{n \to \infty} n^a \beta_1(n) < \infty$, here $(\beta_1(n))_n$ denotes the sequence of strong $\beta_1$-mixing coefficients of $(Z_i)_{i \in \mathbb Z}$
    
        \item $\frac{1}{\sqrt n}\sum_{i=1}^n (|F_\mu(Z_i)|^{1/p}-\E(|F_\mu(Z_i)|^{1/p}))$ ~ does not converge in distribution to a Gaussian law.
    \end{enumerate} 
\end{prp}

Let us see how Proposition \ref{Lp_thm:optimality} can be deduced from Theorem 5 in \cite{DMR1994}. Let $\Tilde F \eqdef \abs{ F_\mu \circ F^{-1}(.)}^{1/p}$. Note that $\int_0^1 u^{-1/a} Q_{Y_{p,\mu}}^2(u) \, du$ is convergent if and only if, with $\Tilde F$ defined above, $\int_0^{1/2} u^{-1/a} \Tilde F^2(u) \, du$ and $\int_0^{1/2} u^{-1/a} \Tilde F^2(1-u) \, du$ are convergent. In other words, in our setting one of the two previous integrals is infinite. Let us assume for instance that $\int_0^{1/2} u^{-1/a} \Tilde F^2(1-u) \, du = +\infty$. Theorem 5 in \cite{DMR1994}, applied to the function $f : x \mapsto \Tilde F(1-x)$, asserts that there exists a stationary Markov chain $(U_i)_{i \in \Z}$ with uniform marginal distributions on $[0,1]$ such that $(\beta_{1,U}(n))_{n \geq 0}$ satisfies \textit{(i)} and 
    $$\frac{1}{\sqrt n}\sum_{i=1}^n (|F_\mu\circ F^{-1}(1-U_i)|^{1/p}-\E(|F_\mu\circ F^{-1}(1-U_i)|^{1/p}))$$
does not converge in distribution to a Gaussian distribution. 
Furthermore, setting $Z_i = F^{-1}(1-U_i)$, the Markov chain $(Z_i)_{i \in \mathbb Z}$ admits F as marginal distribution function and the same mixing coefficients as $(U_i)_{i \in \mathbb Z}$ so that $(i)$ is verified. 
\\

Furthermore, analogously to the first section, we get sufficient conditions for (\ref{Lp_thm:H}) to be satisfied.
\begin{cor}
    Let $(X_i)_{i \in \mathbb Z}$ and $(Y_i)_{i \in \mathbb Z}$ be as defined in (\ref{Lp_condition_F}) and (\ref{Lp_def_X}), and $(\mathcal F _i)_{i \in \Z}$ be a filtration to which $(X_i)_{i \in \Z}$ is adapted. 
    Then under one of these conditions: 
    \begin{enumerate}[label = (\roman*)]
        \item there exists $r > 2$ such that $\E(Y_{p,\mu}^{r}) < \infty$ and $\sum_{n \geq 0} (n+1)^{2/(r-2)} \Tilde \beta_{1,Y}(n) < \infty$

        \item there exists $r > 2$ and $c > 0$ such that for any $x$, $\Proba\bigp{ Y_{p,\mu} > x} \leq \bigp{\frac{c}{x}}^r$ and $\sum_{n \geq 0} \Tilde \beta_{1,Y} (n)^{1-2/r} < \infty$

        \item there exists $a > 0$ and $\tau > 0$ such that $\E \bigcro{Y_{p,\mu}^2 (\log (1+Y_{p,\mu}))^a } < \infty$ and $\Tilde \beta_{1,Y}(n) = O\bigp{e^{-\tau n^{1/a}}}$
    \end{enumerate}
    the condition (\ref{Lp_thm:H}) is verified, so Corollary \ref{Lp_thm} applies. 
\end{cor}
~\\

\noindent\textit{Application to the empirical process in $L^p([0,1],\lambda), p\geq 2$,  for intermittent maps.}
\\
For $\gamma \in ]0,1[$, let $T_\gamma: [0,1] \to [0,1]$ be the intermittent map defined by \cite{LSV1999} as follows: 
\begin{align}
    T_\gamma(x) = \left\{
    \begin{array}{ll}
        x(1+2^\gamma x^\gamma) & if \, x \in [0, 1/2[  \\
        2x - 1 &  if \, x \in [1/2, 1]
    \end{array}
    \right. .
\end{align}
If there is no confusion we write $T$ for the sake of clarity. As shown in \cite{LSV1999}, for all $\gamma \in ]0,1[$, there exists a unique absolutely continuous $T_\gamma$-invariant probability measure $\nu_\gamma$ (or simply $\nu$) on [0,1] whose density $h_\gamma$ satisfies: there exist two finite constants $c_1,c_2 > 0$ such that for all $x \in [0,1]$, $c_1 \leq x^\gamma h_\gamma(x) \leq c_2$. Let us fix $\gamma$ and consider $K$ the Perron-Frobenius operator of $T$ with respect to $\nu$ defined by 
\begin{align}
    \nu(f\circ T.g) = \nu(f.Kg) , \; \textrm{ for any }f,g \in \mathbb L^2(\nu).
\end{align}
Then, by considering $(X_i)_{i \in \mathbb Z}$  a stationary Markov chain with invariant measure $\nu$ and transition kernel $K$, for any positive integer $n$, on the probability space $([0,1], \nu)$, $(T, T^2, \cdots, T^n)$ is distributed as $(X_n, X_{n-1}, \cdots, X_1)$ (see for instance Lemma XI.3 in \cite{HH2001}). Consequently, the two following empirical processes have the same distribution 
\begin{itemize}
    \item $\left\{G_n(t) = \frac{1}{\sqrt n}\sum_{k=1}^n [\indic_{T^k \leq t} - F(t)] \, ; \, t \in [0,1]\right\}$
    \item $\left\{L_n(t) = \frac{1}{\sqrt n}\sum_{k=1}^n [\indic_{X_k \leq t} - F(t)] \, ; \, t \in [0,1]\right\}$
\end{itemize}
where $F(t) = \nu([0,t])$. 
Since $\nu$ is supported on $[0,1]$, condition (\ref{Lp_thm:H}) reads as $\sum_{n \geq 0} \Tilde \beta_{1,X}(n) < \infty$. Now, from \parencite[Proposition 6.2]{DDT2015}, we have the upper bound 
\begin{align}
    \Tilde \beta_{1,X}(n) \leq \frac{C}{(n+1)^{(1-\gamma)/\gamma}}. 
\end{align}
Hence, applying Corollary \ref{Lp_thm}, we derive that for any $\gamma \in ]0,1/2[$ and any $p \geq 2$
\begin{align}
    \{ \, G_n(t) \, : t \in [0,1]\} \xrightarrow[n \to \infty]{\law} \{G(t) \,:\, t \in [0,1]\} \; \textrm{in } L^p([0,1], \lambda)
\end{align}
where $G$ is a Gaussian process. 
\\
We could also consider unbounded but monotonic observables $\varphi$ of the iterates such as $\varphi(x)=1/x^\alpha$ or $\varphi(x)=1/(1-x)^\alpha$ since for such functions, the $\Tilde \beta_1$-coefficients of $\bigp{\varphi(X_i)}_{i \in \Z}$ are of the same order than the initial ones. More precisely taking into account the behavior of the density $h_\gamma$ of $\nu$, one can prove that condition (\ref{Lp_thm:H}) reads as $\alpha < \frac{p}{2}(1-2\gamma)$ when $\varphi(x)=1/x^\alpha$ and $\alpha < \frac{p}{2}\frac{1-2\gamma}{1-\gamma}$ when $\varphi(x)=1/(1-x)^\alpha$. 

\section{Proofs}

We start this section by a general CLT for random variables taking values in 2-smooth Banach spaces. It will be a building block in the proof of Theorem \ref{thm_ppal} and has interest in itself. 

\subsection{A general result}

\begin{thm}\label{thm_ros}
    Let $\mathbb B$ be a 2-smooth Banach space and let $(X_i)_{i \in \Z}$ be a stationary sequence of $\mathbb B$-valued centered random variables, adapted to a non-decreasing and stationary filtration $(\mathcal F _i)_{i \in \Z}$, and such that $\E(\norm{X_0}_{\mathbb B}^2) < \infty$. 
    Set $S_n \eqdef \sum\limits_{k = 1}^n X_k$ and assume that  
    \begin{align}
        &\textrm{for any } x^* \in \mathbb B^*, \, \bigp{\frac{[x ^* (S_n)]^2}{n}}_n \textrm{is an uniformly integrable family} \label{thm_ros:H1},
        \\
        &\frac{1}{n} \E \bigcro{ \norm{\E(S_n | \calF _0)}_{\mathbb B}^2 } \xrightarrow[n \to +\infty]{} 0 \label{thm_ros:H2},
        \\
        &\textrm{for any } x^* \in \mathbb B ^*, \textrm{ there exists } \sigma^2(x^*) \textrm{ such that } \E \abs{\E\bigp{\left.\frac{[x ^* (S_n)]^2}{n} \right\lvert \calF _0} - \sigma^2(x^*)} \xrightarrow[n \to +\infty]{} 0 \label{thm_ros:H3},
        \\
        \textrm{and there}&\textrm{ exists } (F_l)_l \textrm{ a sequence of finite dimensional subspaces of } \mathbb B \textrm{ such that } \nonumber
        \\
        &\limsup\limits_{n \to +\infty} \frac{1}{n} \E [q_{F_l}^2(S_n)] \xrightarrow[l \to +\infty]{} 0 \label{thm_ros:H4},
        \\
        \textrm{where } q_F&(x)  \eqdef \inf\{ \norm{x-y}_{\mathbb B} : y \in F\} \textrm{ for any }~ x \in \mathbb B
        .
        \nonumber
    \end{align}
    Then $\bigp{\frac{1}{\sqrt n}S_n}_{n \geq 1}$ converges in distribution to $G$, where $G$  is a $\mathbb B$-valued Gaussian random variable whose law $\gamma$ is such that for any $x^* \in \mathbb B ^*$, $\Hat{\gamma}(x^*) = e^{-\sigma^2(x^*)/2}$. The covariance operator of $G$ is then given by: $K_G(x^*, y^*)= \frac{1}{2}\bigp{\sigma^2(x^*+y^*) - \sigma^2(x^*)-\sigma(y^*)}$ for any $x^*, y^* \in \mathbb B^*$.
\end{thm}
~\\
\begin{rem}\label{thm_ros_cond_equiv}
    Conditions (\ref{thm_ros:H1}) and (\ref{thm_ros:H2}) can be replaced by the two following conditions: 
        \begin{flalign}
        \qquad \qquad &\bigp{\frac{\norm{S_n}_{\mathbb B}^2}{n}}_n \textrm{is an uniformly integrable family} \label{thm_ros:H1b}, &&
        \\
        &\frac{1}{\sqrt n} \E \norm{\E(S_n | \calF _0)}_{\mathbb B} \xrightarrow[n \to +\infty]{} 0
        . 
        \label{thm_ros:H2b} &&
        \end{flalign}
\end{rem}
~\\
\textit{Proof of Remark \ref{thm_ros_cond_equiv}. }
    Let $x^*$ be an element of $\mathbb B ^*$. Then, there exists $c(x^*) > 0$ such that for any $x \in \mathbb B$, $\abs{x^*(x)}\leq c(x^*)\norm{x}_{\mathbb B}$. Hence, 
    \begin{align*}
        \forall n \in \N ^*, \, \abs{\frac{[x ^* (S_n)]^2}{n}} \leq c(x^*)^2\, \frac{\norm{S_n}_{\mathbb B}^2}{n}
        .
    \end{align*}
    So, condition (\ref{thm_ros:H1b}) implies that $\bigp{\frac{[x ^* (S_n)]^2}{n}}_n$ is an uniformly integrable family. 
    \\
    On another hand, for any $A>0$, 
    \begin{align*}
        \frac 1 n \norm{\E(S_n | \calF _0)}_{\mathbb B}^2
        & \leq 
        \bigp{\frac{1}{\sqrt n} \norm{\E(S_n | \calF _0)}_{\mathbb B}\indic_{\norm{\condEp{\frac{S_n}{\sqrt n}}{\calF _0}}_{\mathbb B} \leq A}} ^2  + \bigp{\frac{1}{\sqrt n} \norm{\E(S_n | \calF _0)}_{\mathbb B} \indic_{\norm{\condEp{\frac{S_n}{\sqrt n}}{\calF _0}}_{\mathbb B} > A} }^2
        .
    \end{align*}
    Clearly,
    \begin{align*}
        \bigp{\frac{1}{\sqrt n} \norm{\E(S_n | \calF _0)}_{\mathbb B}\indic_{\norm{\condEp{\frac{S_n}{\sqrt n}}{\calF _0}}_{\mathbb B} \leq A}} ^2
        \leq \frac{A}{\sqrt n}\norm{\E(S_n | \calF _0)}_{\mathbb B}
        ,
    \end{align*}
    which converges to 0 in $\mathbb L^1$ by condition (\ref{thm_ros:H2b}). 
    On another hand, by  Lemma 6.3 in \cite[]{DMP2014}, we have for any $\mathbb B$-valued r.v. $Y$, any $\eps >0$ and any $\sigma$-algebra $\mathcal B$, 
    \begin{align}\label{lemma_contraction_indic_Econd}
        \E \bigp{
            \norm{Y}_{\mathbb B}^2 \indic_{\norm{\condEp{Y}{\mathcal B}}_{\mathbb B} > 2 \varepsilon}
        }
        \leq 
        \E \bigp{
            \norm{Y}_{\mathbb B}^2 \indic_{\norm{Y}_{\mathbb B} > \varepsilon}
        } 
        .
    \end{align}
    \\
    Hence, applying Jensen's inequality and taking into account (\ref{lemma_contraction_indic_Econd}), we get 
    \begin{align*}
        \E \bigcro{\bigp{\frac{1}{\sqrt n} \norm{\E(S_n | \calF _0)}_{\mathbb B} \indic_{\norm{\condEp{\frac{S_n}{\sqrt n}}{\calF _0}}_{\mathbb B} > A} }^2}
        &\leq 
        \E \bigp{\frac{\norm{S_n}_{\mathbb B} ^2}{n} \indic_{\norm{\condEp{\frac{S_n}{\sqrt n}}{\calF _0}}_{\mathbb B} > A} }
        \\
        & \leq 
        \E \bigp{\frac{\norm{S_n}_{\mathbb B} ^2}{n} \indic_{\norm{S_n}_{\mathbb B}/\sqrt n > A/2}}
        ,
    \end{align*}
    which converges to zero by (\ref{thm_ros:H1b}) by first letting n tend to infinity and after A.
    So, overall, (\ref{thm_ros:H1b}) together with (\ref{thm_ros:H2b}) imply (\ref{thm_ros:H2}). 
\finpreuve

\subsection{Proof of Theorem \ref{thm_ros}}
   To prove it, we follow the method of proof given in the proof of \parencite[Theorem 4.31]{MPU2019}. This method consists in constructing big blocks of random variables and then in approximating them by a triangular array of martingales to which a CLT is applied. In the Banach space setting, to prove Theorem \ref{thm_ros} we shall rather use \parencite[Theorem 5]{Ros1982} instead of \parencite[Theorem 2.4]{MPU2019}. 
\\ 
    For any fixed positive integer $m$, set $p = \lfloor n/m\rfloor$ and let  
    \begin{align*}
        X_{n,j} \eqdef \sum\limits_{k = p(j-1)+1}^{pj} \frac{1}{\sqrt n} X_k 
        \quad \textrm{and} \quad 
        \calF_{n,j} \eqdef \calF _{jp}.
    \end{align*}
    Let us consider the conditionally centered random variables
    \begin{align*}
        \Tilde{X_{n,j}} \eqdef X_{n,j} - \condEp{X_{n,j}}{\calF_{n,j-1}},
    \end{align*}
    so that $(\Tilde{X_{n,j}})_{1\leq j \leq n}$ is a triangular array of martingale differences adapted to the array of filtrations $(\calF_{n,j})_{1\leq j \leq n}$.
    \\
    As done in \cite[proof of Theorem 4.31, p. 133]{MPU2019}, let consider $(m_n)_n$ a sequence sufficiently slow growing such that, as $n$ tends to infinity   
    \begin{align}\label{thm_ros_proof:mn_H1}
        m = m_n \xrightarrow[]{} + \infty ,
        \qquad
        \frac{m}{\sqrt n} \xrightarrow[]{} 0
    \end{align}
    and,
    \begin{align}\label{thm_ros_proof:mn_H2}
        \sqrt{m} \, \condEp{\frac{1}{\sqrt p}S_p}{\calF_0} \xrightarrow[n \to +\infty]{} 0 \quad in \; \mathbb L^1(\mathbb B)
        .
    \end{align}
    We have in particular, from (\ref{thm_ros_proof:mn_H1}), that $\frac{mp}{n} \xrightarrow[]{} 1$. 
    \\
    Note that 
    \begin{align*}
        \sum\limits_{j=1}^m \Tilde{X_{n,j}}
        &= \sum\limits_{j = 1}^m X_{n,j} -  \sum\limits_{j = 1}^m \condEp{X_{n,j}}{\calF_{n,j-1}}
        \\ 
        &= \frac{1}{\sqrt n} S_{n} + \frac{1}{\sqrt n}(S_{pm}-S_n) - \sum\limits_{j = 1}^m \condEp{X_{n,j}}{\calF_{n,j-1}}.
    \end{align*}
    On the one hand, 
    \begin{align*}
        \E \norm{\frac{1}{\sqrt n}(S_{pm}-S_n)}_{\mathbb B} 
        &\leq \frac{m}{\sqrt n} \, \E \norm{X_0}_{\mathbb B},
    \end{align*}
    which converges to 0 as $n$ tends to $+\infty$ from (\ref{thm_ros_proof:mn_H1}). 
    On another hand, bearing in mind stationarity,
    \begin{align*}
        \E \norm{\sum\limits_{j = 1}^m \condEp{X_{n,j}}{\calF_{n,j-1}}}_{\mathbb B} \leq \sqrt{\frac{mp}{n}}\, \sqrt{m} \, \E \norm{\condEp{\frac{S_p}{\sqrt p}}{\calF_0}}_{\mathbb B},
    \end{align*}
    which converges to 0 as $n$ tends to $+\infty$ combining conditions (\ref{thm_ros_proof:mn_H1}) and (\ref{thm_ros_proof:mn_H2}). 
    \\
    So, overall, Theorem \ref{thm_ros} will be proved if one can show that $\sum\limits_{j = 1}^m \Tilde{X_{n,j}}$ converges in distribution to $G$. 
    To this end, we shall use the following result due to \parencite[Theorem 5]{Ros1982}. 
    \begin{lem}[Rosinski]
        Let $\mathbb B$ be a $2$-smooth Banach space. Let $(Y_{n,j})_{n \geq 1, j \leq k_n}$ be a martingale difference array associated to $\sigma$-fields $(\mathcal{A}_{n,j})_{n \geq 1, j \leq k_n}$ such that $\E\abs{x^*(Y_{n,j})}^2 < \infty$ for any $x^* \in \mathbb B^*$. Assume that
        \begin{enumerate}[label = (\textit{\roman*)}]
            \item\label{ros:H1} there exists $\psi : \mathbb B^* \to \R_+$ such that for any $x^* \in \mathbb B^*$, $\sigma_n^2(x^*) \xrightarrow[]{\Proba} \psi(x^*)$ as $n \to +\infty$, where $\sigma_n^2(x^*) = \sum_{j =1}^{k_n} \condEp{[x^*(Y_{n,j})]^2}{\mathcal{A}_{n,j-1}}$
            
            \item\label{ros:H2} for any $x^* \in \mathbb B^*$ and any $\eps >0$, $\sum_{j = 1}^{k_n} \condEc{x^*(Y_{n,j})^2 \indic_{\abs{x^*(Y_{n,j})} > \eps}}{\mathcal{A}_{n, j-1}} \xrightarrow{\Proba} 0 $ as $n \to +\infty$
            
            \item\label{ros:H3} there exists a sequence $(F_l)_l$ of finite dimensional subspaces of $\mathbb B$ such that 
                \begin{align*}
                    \limsup_{n \to +\infty} \sum_{j = 1}^{k_n} \condEc{q_{F_l}^2(Y_{n,j})}{\mathcal{A}_{n, j-1}} \xrightarrow{\Proba} 0 \textrm{ as } l \to +\infty.
                \end{align*}
        \end{enumerate}
        Then $\bigp{\sum_{j=1}^{k_n}Y_{n,j}}_{n \geq 1}$ converges in distribution to a centered Gaussian r.v. whose distribution $\gamma$ verifies $\hat \gamma(x^*) = e^{-\psi(x^*)}$ for any $x^* \in \mathbb B ^*$.
    \end{lem}
    Let us apply it to $(\Tilde{X_{n,j}})_{n\geq 1, j\leq m_n}$ and $(\mathcal{F}_{n,j})_{n \geq 1, j \leq m_n}$.
    First, since for any integer $k$, $\E(\norm{X_k}_{\mathbb B}^2)$ is finite, it appears immediately that for any $x^* \in \mathbb B ^*$ and any $n,j$, $\E [x^* (\Tilde{X_{n,j}})]^2 < \infty$. 
    \\
    
    Let us begin with the proof of condition \ref{ros:H1}. 
    Let $x^*$ be a linear form defined on $\mathbb B$. 
    Using the same notations as in Rosinski's paper, we write 
    \begin{align*}
        \sigma_n^2(x^*) \eqdef \sum\limits_{j = 1}^m \condEc{[x^* (\Tilde{X_{n,j}})]^2}{\calF_{n,j-1}}
        = \sum\limits_{j = 1}^m \bigp{\condEc{[x^* (X_{n,j})]^2}{\calF_{n,j-1}} - \condEc{x^* (X_{n,j})}{\calF_{n,j-1}}^2}
        .
    \end{align*}
    Hence, by stationary and condition (\ref{thm_ros_proof:mn_H1}), we derive 
    \begin{align}\label{thm_ros_proof:sigman_decomposition}
        \limsup\limits_{n \to +\infty} \E \,\abs{\sigma _n^2 (x^*) - \sigma ^2(x^*)}
        \leq \limsup\limits_{p \to + \infty} \E\abs{\condEp{\frac{[x^* ( S_p)]^2}{p}}{\calF_0} - \sigma ^2(x^*)} + \limsup\limits_{p \to +\infty} \E \bigp{ \abs{\condEp{\frac{x^* (S_p)}{\sqrt p}}{\calF_0}} ^2 }
        .
    \end{align}
    Then, combining (\ref{thm_ros_proof:sigman_decomposition}), (\ref{thm_ros:H2}) and (\ref{thm_ros:H3}), we get
    \begin{align*}
        \sigma_n^2(x^*) \xrightarrow[n \to +\infty]{\mathbb L ^1} \sigma^2 (x^*) .
    \end{align*} 
    This ends the proof of \ref{ros:H1}.
    \\
    
    Let us now check condition \ref{ros:H2}. 
    Let $\eps$ be a positive real number. For any $n \in \N ^*$, in accordance with \cite[Lemma~2.28]{MPU2019} and by stationarity  
    \begin{align*}
        &\E \abs{\sum\limits_{j = 1}^m \condEp{[x^*(\Tilde{X_{n,j}})]^2 \, \indic_{\abs{x^*(\Tilde{X_{n,j}})}>\eps}}{\calF_{n,j-1}}}
        \\
        & \hspace{20pt}
        \leq 12 \E \, \sum\limits_{j = 1}^m \condEp{[x^*( X_{n,j} )]^2 \, \indic_{\abs{x^*( X_{n,j} ) }>\eps/4}}{ \calF_{n,j-1} }
        = 12\frac{mp}{n} \E \bigcro{\frac{[x^*( S_p )]^2}{p} \, \indic_{\abs{\frac{x^*( S_p )}{\sqrt p} }>\frac{\sqrt n}{4\sqrt p}\eps}}
        .
    \end{align*}
    Then, from (\ref{thm_ros:H1}) and since $\frac{mp}{n} \to 1$ as $n$ tends to $+\infty$, 
    \begin{align*}
        \sum\limits_{j = 1}^m \condEp{[x^*(\Tilde{X_{n,j}})]^2 \, \indic_{\abs{x^*(\Tilde{X_{n,j}})}>\eps}}{\calF_{n,j-1}} \xrightarrow[n \to +\infty]{} 0 \quad in \; \mathbb L ^1
        . 
    \end{align*}
    ~\\
    
    It remains to verify condition \ref{ros:H3}.
    It is easy to see that for any $x,y \in \mathbb B$, $q_{F_l}(x+y) \leq q_{F_l}(x) + \norm{y}_{\mathbb B}$. Therefore
    \begin{align*}
        q_{F_l}(\Tilde{X_{n,j}}) \leq q_{F_l}(X_{n,j}) + \norm{\condEp{X_{n,j}}{\calF_{n,j-1}}}_{\mathbb B}
        ,
    \end{align*}
    implying that
    \begin{align*}
        q_{F_l}^2(\Tilde{X_{n,j}}) \leq 2 q_{F_l}^2(X_{n,j}) + 2\norm{\condEp{X_{n,j}}{\calF_{n,j-1}}}_{\mathbb B}^2
        .
    \end{align*}
    On the one hand, taking into account stationarity, we derive 
    \begin{align*}
        \E \, \sum\limits_{j = 1}^m \norm{ \condEc{ X_{n,j} }{\calF_{n,j-1}} }_{\mathbb B}^2
        &= \frac{mp}{n} \E \norm{\condEp{\frac{S_p}{\sqrt p}}{\calF_0}}_{\mathbb B}^2 
        ,
    \end{align*}
    which converges to 0 as $n$ tends to infinity by conditions (\ref{thm_ros:H2}) and (\ref{thm_ros_proof:mn_H1}). 
    \\
    On the other hand, by stationarity  
    \begin{align*}
         \sum\limits_{j = 1}^m  \E \bigcro{ q_{F_l}^2(X_{n,j}) } = m \E \bigcro{ q_{F_l}^2(X_{n,1}) } = \frac{mp}{n}\frac{1}{p} \E[q_{F_l}^2(S_p)]
         .
    \end{align*}
    Therefore, combining (\ref{thm_ros:H4}) and (\ref{thm_ros_proof:mn_H1}), we derive  
    \begin{align*}
        \lim\limits_{l \to +\infty} \, \limsup\limits_{n \to +\infty} \sum\limits_{j = 1}^m  \E \bigcro{ q_{F_l}^2(X_{n,j}) } = 0
        .
    \end{align*}
    Hence, $(\Tilde{X_{n,j}})_{n,j}$ satisfies the conditions of \parencite[Theorem 5]{Ros1982}. 
    This ends the proof of Theorem \ref{thm_ros}.

\subsection{Proof of Theorem \ref{thm_ppal}}

~ 
    \\ 
    We shall apply Theorem \ref{thm_ros}.
    \\
    ~
    Let start by showing that (\ref{thm_ppal:H}) implies (\ref{thm_ros:H1}). 
    Denote for any nonnegative integer $k$,
    \begin{align}\label{thm_ppal_proof:def_Y}
        Y_k = x^* (X_k)
        \qquad
        \textrm{ and } \qquad
        T_k = \sum\limits_{i = 1}^k Y_i
        .
    \end{align}
    Condition (\ref{thm_ppal:H}) implies that for any $x^* \in \mathbb B ^*$, $(x^*(X_0)\condEp{x^*(S_n)}{\calF_0})_n$ converges in $\mathbb L^1$.  Hence, applying \parencite[Proposition 1(b)]{DR2000}, it follows that $\bigp{\frac{T_n^2}{n}}_n$ is an uniformly integrable family, which proves (\ref{thm_ros:H1}).
    \\
    
    We prove now that under (\ref{thm_ppal:H}), (\ref{thm_ros:H4}) is verified. 
    Let $(e_n)_{n \in \N}$ be a Schauder basis of $\mathbb B$. For any $l$, denote by $P_l$ the projection on the subspace $F_l$ generated by the $l$ first vectors of the basis. 
    For any integers $n$ and $l$, since $P_l(S_n) \in F_l$, 
    \begin{align*}
        \frac{1}{n} \E \, q_{F_l}^2 (S_n) \leq \frac{1}{n}\E \norm{(\textrm{id} - P_l)S_n}_{\mathbb B}^2
        .
    \end{align*}
        From Theorem 2.1 together with Lemma 1.1 and Remark 2.1 in \cite{DM2015} with $p = 2$, we get that there exists $c>0$ such that for any positive integer $n$  
        \begin{align*}
            \E \norm{(\textrm{id} - P_l)S_n}_{\mathbb B}^2 \leq c \sum\limits_{i = 1}^n \max\limits_{i \leq j \leq n} \E \bigp{\norm{(\textrm{id} - P_l) X_i}_{\mathbb B} \norm{\sum\limits_{k=i}^j \condEc{(\textrm{id} - P_l) X_k}{\calF _i}}_{\mathbb B}}
            .
        \end{align*}
        Hence, by using stationarity, we derive 
        \begin{align*}
            \E \norm{(\textrm{id} - P_l)S_n}_{\mathbb B}^2 
            &
            \leq c \sum\limits_{i = 1}^n \max\limits_{i \leq j \leq n} \E \bigp{\norm{(\textrm{id} - P_l) X_0}_{\mathbb B} \norm{\sum\limits_{k=0}^{j-i} \condEc{(\textrm{id} - P_l) X_k}{\calF_0}}_{\mathbb B}} \nonumber
            \\
            & 
            \leq 
            nc \max\limits_{1 \leq j \leq n} \E \bigp{\norm{(\textrm{id} - P_l) X_0}_{\mathbb B} \norm{\sum\limits_{k=0}^{j-1} \condEc{(\textrm{id} - P_l) X_k}{\calF_0}}_{\mathbb B}}\nonumber
            .
    \end{align*}
    Furthermore, $(P_l)_l$ is uniformly bounded for the operator norm, then so is $(\textrm{id} - P_l)_l$. Hence there exists $C >0$ such that for any $l$ and any $n > N$,
    \begin{align}
            &
            \frac{1}{n} \E \norm{(\textrm{id} - P_l)S_n}_{\mathbb B}^2  \nonumber
            \\
            & \hspace{1cm}
            \leq 
            C \, \E \bigp{\norm{(\textrm{id} - P_l) X_0}_{\mathbb B} \norm{\sum\limits_{k=0}^{N} \condEp{X_k}{\calF_0}}_{\mathbb B}}
            + C \max\limits_{N+1 \leq j \leq n}
            \E \bigp{\norm{X_0}_{\mathbb B}  \norm{\sum\limits_{k=N+1}^{j-1} \condEp{X_k}{\calF_0}}_{\mathbb B}}
            \label{thm_ppal_proof:majoration_norm(id-Pl)}
            .
        \end{align}
    \noindent
    The second term in the right hand side converges to zero by condition (\ref{thm_ppal:H}) by letting first $n$ tend to $+\infty$ and after $N$. 
    \\
    Since $\E(\norm{X_0}_{\mathbb B}^2) < \infty$ and $\norm{(\textrm{id}-P_l)X_0}_{\mathbb B}$ converges a.s. to zero as $l$ tends to $+\infty$, the first term in the right hand side of (\ref{thm_ppal_proof:majoration_norm(id-Pl)}) converges to zero by letting first $l$ tend to $+\infty$. So, overall, 
    \begin{align}\label{thm_ppal_proof:convergence_norm(id-Pl)}
        \lim_{l \to +\infty} \limsup_{n \to +\infty} \frac{1}{n} \E (\norm{(\textrm{id}-P_l)S_n}_{\mathbb B}^2) = 0,
    \end{align}
    and (\ref{thm_ros:H4}) follows.
    \\
    
    We prove now that (\ref{thm_ppal:H}) implies (\ref{thm_ros:H2}). 
    For any positive integers $n$ and $l$, we have  
    \begin{align*}
        \frac{1}{n}\E \bigcro{ \norm{\condEp{S_n}{\calF_0}}_{\mathbb B}^2 }
        &\leq \frac{2}{n} \E \bigcro{ \norm{\condEp{ P_l \, S_n}{\calF_0}}_{\mathbb B}^2 } + \frac{2}{n} \E \bigcro{ \norm{\condEc{(\textrm{id}-P_l) \, S_n}{\calF_0}}_{\mathbb B}^2 }
        .
    \end{align*}
    The second term in the right hand side is going to zero as $n$ tends to infinity by using (\ref{thm_ppal_proof:convergence_norm(id-Pl)}) and Jensen's inequality. 
    \\
    Furthermore, to prove the convergence of $\frac{2}{n} \E \bigcro{ \norm{\condEp{ P_l \, S_n}{\calF_0}}_{\mathbb B}^2 }$ to 0 as $n$ tends to $+\infty$, it is sufficient to prove that for any $x^* \in \mathbb B ^*$ we have  $\frac{1}{n} \E \bigcro{\condEp{ x^* S_n }{\calF_0}^2} \rightarrow 0$.
    Since $\bigp{\frac{[x ^* (S_n)]^2}{n}}_n$ is uniformly integrable, it is then enough to prove that  
    \begin{align*}
        \frac{1}{\sqrt n}\E \, \abs{ \condEp{x^*(S_n)}{\calF_0} } \xrightarrow[n \to +\infty]{} 0
        .
    \end{align*}
    This holds under the condition $\bigp{x^*(X_0) \, x^*\bigp{\condEp{S_n}{\calF_0}}}_n$ converges in $\mathbb L^1$ by the arguments developed in Step 3 of the proof of Theorem 4.18 in \cite{MPU2019}. 
    \\
    
    It remains to check Condition (\ref{thm_ros:H3}). 
    Applying $x^*$, we place ourselves in the well-known context of real random variables. From Condition (\ref{thm_ppal:H}) and using the same notations as in (\ref{thm_ppal_proof:def_Y}), $(Y_0\condEp{T_n}{\calF_0})_n$ converges in $\mathbb L^1$ then the series of covariance associated to $(Y_k)_k$ converges. Let us denote it $\sigma^2(x^*)$. \\
    From \parencite[Proof of Theorem 4.18, Step 4]{MPU2019}, since the sequence $(X_n)_{n \in \mathbb Z}$ is ergodic, we get  
    \begin{align*}
        \frac{1}{n} \condEp{T_n ^2}{\calF_0} \xrightarrow[n \to +\infty]{\mathbb L ^1} \sigma^2(x^*)
        .
    \end{align*}
    That is  
    \begin{align*}
        \condEp{\frac{[x ^* (S_n)]^2}{n}}{\calF_0} \xrightarrow[n \to +\infty]{\mathbb L ^1} \sigma^2(x^*) \eqdef \sum\limits_{k \in \Z} \textrm{Cov}(x^*(X_0), x^* (X_k))
        .
    \end{align*}
    This ends the proof of Theorem \ref{thm_ppal}. 

\subsection{Proof of Corollary \ref{Lp_thm}}

Note that when $\mu$ is finite, the result follows immediately. Assume from now that $\mu$ is not finite. 
\\
Let start by noting that 
        $\norm{X_0}_{p, \mu} \leq Y_{p, \mu} + \E Y_{p,\mu}$, so that 
        \begin{align}\label{Lp_thm_proof:link_normX_Ypmu}
            Q_{\norm{X_0}_{p, \mu}} \leq Q_{Y_{p,\mu}} + \E Y_{p, \mu}
            \quad\textrm{ and }\quad
            G_{\norm{X_0}_{p, \mu}}(.) \geq G_{Y_{p,\mu}}(./2),
        \end{align}
where $\norm{X_0}_{p, \mu} = \bigp{\int_{\R} \abs{X_0(t)}^p \, d\mu(t)}^{1/p}$. 
Indeed, by definition of $X_0$ and as an application of Minkowski's inequality, 
    \begin{align*}
        \norm{X_0}_{p,\mu} 
        &= 
        \bigp{ \int_{-\infty}^0 \abs{\indic_{Y_0 \leq t} - F(t)}^p \, d\mu(t) + \int^{+\infty}_0 \abs{\indic_{Y_0 > t} +1 - F(t)}^p \, d\mu(t) }^{1/p}
        \\
        &\leq 
        \bigp{ \int_{-\infty}^0 \indic_{Y_0 \leq t} \, d\mu(t) + \int^{+\infty}_0 \indic_{Y_0 > t} \, d\mu(t) }^{1/p} + \bigp{ \int_{-\infty}^0 \abs{F(t)}^p \, d\mu(t) + \int^{+\infty}_0 \abs{1 - F(t)}^p \, d\mu(t) }^{1/p}
        \\
        &\leq
        Y_{p,\mu} + \E Y_{p,\mu}
        .
    \end{align*}
    \noindent
    The inequality for the quantile function follows immediately. Let us prove the last inequality.  
    \\
    For any $x \in [0,1]$,  
    \begin{align*}
        \int_0^x Q_{\norm{X_0}_{p,\mu}}(u) \, du
        \leq \int_0^x Q_{Y_{p,\mu}}(u) \, du + x \, \E Y_{p,\mu}
        \leq \int_0^x Q_{Y_{p,\mu}}(u) \, du + x \, \int_0^1 Q_{Y_{p,\mu}}(u) \, du
        .
    \end{align*}
    Hence, as $Q_{Y_{p,\mu}}$ is non-increasing, 
        $\int_0^x Q_{\norm{X_0}_{p,\mu}}(u) \, du \leq 2\int_0^{x/2} Q_{Y_{p,\mu}}(u)$. Thus, $G_{\norm{X_0}_{p, \mu}}(.) \geq G_{Y_{p,\mu}}(./2)$ and (\ref{Lp_thm_proof:link_normX_Ypmu}) is established. 
\\
    From (\ref{Lp_thm_proof:link_normX_Ypmu}) and after a change of variables, for any $a>0$,
    \begin{align}\label{Lp_thm_proof:ineq_ppal_step1}
        \int_0^{\E\norm{\condEp{X_n}{\calF_0}}_{p,\mu}} Q_{\norm{X_0}_{p,\mu}}& \circ G_{\norm{X_0}_{p,\mu}} (u) \, du \nonumber
         \\ 
          &\leq 
        2\int_0^{\E\norm{\condEp{X_n}{\calF_0}}_{p,\mu}/2} Q_{Y_{p,\mu}} \circ G_{Y_{p,\mu}} (u) \, du + \E\norm{\condEp{X_n}{\calF_0}}_{p,\mu}.\E Y_{p,\mu}
        \nonumber
        \\
        &\leq 
        2a\int_0^{G_{Y_{p,\mu}}(\E\norm{\condEp{X_n}{\calF_0}}_{p,\mu}/a)} Q_{Y_{p,\mu}}^2 (u) \, du + \E\norm{\condEp{X_n}{\calF_0}}_{p,\mu}.\E Y_{p,\mu}
        .
    \end{align}
    Let $U$ and $V$ be real r.v.'s such that $U$ is distributed as $\norm{\condEp{X_n}{\calF_0}}_{p,\mu}$, $V$ is distributed as $Y_{p,\mu}$ and $U$ and $V$ are independent. 
    It follows that  
    $\E\norm{\condEp{X_n}{\calF_0}}_{p,\mu}.\E Y_{p,\mu} = \E(UV)$. 
    \\
    Following the proof of Proposition 1, (4.1) in \cite{DD2003} and after taking into account (\ref{Lp_thm_proof:link_normX_Ypmu}), we derive 
    \begin{align}\label{Lp_thm_proof:ineq_ppal_step1_terme2}
        \E(UV) 
        \leq \int_0^{\E\norm{\condEp{X_n}{\calF_0}}_{p,\mu}} Q_{Y_{p,\mu}} \circ G_{\norm{X_0}_{\mathbb B}} (u) \, du 
        \leq 2a\int_0^{G_{Y_{p,\mu}}(\E\norm{\condEp{X_n}{\calF_0}}_{p,\mu}/a)} Q_{Y_{p,\mu}}^2 (u) \, du
        .
    \end{align}
    Finally, (\ref{Lp_thm_proof:ineq_ppal_step1}) together with (\ref{Lp_thm_proof:ineq_ppal_step1_terme2}) imply that 
    \begin{align}\label{Lp_thm_proof:ineq_ppal}
        \int_0^{\E\norm{\condEp{X_n}{\calF_0}}_{p,\mu}} Q_{\norm{X_0}_{p,\mu}}& \circ G_{\norm{X_0}_{p,\mu}} (u) \, du
        &\leq 4a\int_0^{G_{Y_{p,\mu}}(\E\norm{\condEp{X_n}{\calF_0}}_{p,\mu}/a)} Q_{Y_{p,\mu}}^2 (u) \, du .
    \end{align}
    \\

\noindent In the following, we will assume without loss of generality that $\Proba(F_\mu(Y_0)> 0) > 0$ and $\Proba(F_\mu(Y_0)< 0) > 0$, as if it is not the case, in the following calculations some terms disappear making it simpler. 
Let $x > 0$ and $y < 0$, we can write
\begin{align}\label{Lp_thm_proof:case1_majoration_norm(Econd)}
        \norm{\condEp{X_n}{\calF_0}}_{p,\mu}
        &= 
        \left(\int_{-\infty}^y \abs{\condEp{\indic_{Y_n \leq t}}{\calF_0}-F(t)}^p \, d\mu(t) + \int_y^x \abs{\condEp{\indic_{Y_n \leq t}}{\calF_0}-F(t)}^p \, d\mu (t) \right.
        \nonumber 
        \\
        & \qquad \qquad \qquad \qquad \qquad \qquad \qquad \qquad 
        \left. + \int_x^{+\infty} \abs{\condEp{\indic_{Y_n \leq t}}{\calF_0}-1+1-F(t)}^p \, d\mu (t)\right)^{1/p}
        \nonumber 
        \\
        &\leq 
        \bigp{ 
        \int_{-\infty}^y \abs{\condEp{\indic_{Y_n \leq t}}{\calF_0}}^p \, d\mu(t)}^{1/p} + \bigp{\int_{-\infty}^y \abs{\E(\indic_{Y_0 \leq t})}^p \, d\mu(t)
        }^{1/p}
        + F_\mu(x)^{1/p}b_n
        \nonumber
        \\
        + (-F_\mu & (y))^{1/p} b_n
        + \bigp{
        \int^{+\infty}_x \abs{\condEp{\indic_{Y_n > t}}{\calF_0}}^p \, d\mu(t)}^{1/p} + \bigp{ \int^{+\infty}_x \abs{\E(\indic_{Y_0 > t})}^p \, d\mu(t)
        }^{1/p}
        ,
    \end{align}
where $b_n$ is defined in Definition \ref{def_weak_beta}. Considering for any $f$, $\norm{f}_{p, I,\mu} = \bigp{\int_I \abs{f(t)}^p \, d\mu(t)}^{1/p}$ and using Jensen's inequality for $\norm{\cdot}_{p, [x, +\infty[,\mu}$, we get 
\begin{align}
    &\bigp{
        \int^{+\infty}_x \abs{\condEp{\indic_{Y_n > t}}{\calF_0}}^p \, d\mu(t)}^{1/p}
        = \norm{\condEp{\indic_{Y_n > \cdot}}{\calF_0}}_{p, [x, +\infty[,\mu} \nonumber
        \\
        &\hspace{2cm}
        \leq \condEp{\norm{\indic_{Y_n > \cdot}}_{p, [x, +\infty[,\mu}}{\calF_0}
        = \condEc{\bigp{\int^{+\infty}_x \indic_{Y_n > t}\, d\mu(t)}^{1/p}}{\calF_0} 
        \nonumber
    \\
    \textrm{and}&
        \nonumber
    \\
    &\bigp{
        \int^{+\infty}_x \abs{\E(\indic_{Y_0 > t})}^p \, d\mu(t)}^{1/p}
        = \norm{\E\bigp{\indic_{Y_0 > \cdot}}}_{p, [x, +\infty[,\mu} \nonumber
        \\
        &\hspace{2cm}
        \leq \E{\norm{\indic_{Y_0 > \cdot}}_{p, [x, +\infty[,\mu}}
        \leq 
        \E \bigcro{ \bigp{\int^{+\infty}_x \indic_{Y_0 > t}\, d\mu(t)}^{1/p}}.
        \nonumber
\end{align}
We proceed in a similar way with $\norm{\cdot}_{p, ]-\infty, y],\mu}$, so that taking the expectation in (\ref{Lp_thm_proof:case1_majoration_norm(Econd)}), we finally get 
\begin{align}\label{Lp_thm_proof:case1_majoration_E(norm(Econd))}
    &\E \norm{\condEp{X_n}{\calF_0}}_{p,\mu} 
    \\
        &\quad 
        \leq 2\E \bigcro{ \bigp{\int^{+\infty}_x \indic_{Y_0 > t}\, d\mu(t)}^{1/p}}
         + \Tilde \beta_{1,Y}(n) [F_\mu(x)^{1/p} + (-F_\mu(y))^{1/p}]
         + 2 \E \bigcro{\bigp{\int^{y}_{-\infty} \indic_{Y_0 \leq t} \, d\mu(t)}^{1/p}}
         .
    \nonumber
\end{align}

\noindent Note now that
    \begin{align}
        &\E \bigcro{\bigp{ \int_x^{+\infty} \indic_{Y_0 > t} \, d\mu(t)}^{1/p}}
         =
        \E \bigcro{ ((F_\mu(Y_0))_+ - F_\mu(x))_+^{1/p} }
        \nonumber
        \\
        & \qquad \leq 
        \E \bigcro{(F_\mu(Y_0))_+^{1/p} \indic_{(F_\mu(Y_0))_+^{1/p} > F_\mu(x)^{1/p}}}
        \leq 
        \int_0^{Q_{(F_\mu(Y_0))_+^{1/p}}^{-1}(F_\mu (x)^{1/p})} Q_{(F_\mu(Y_0))_+^{1/p}} (u) \, du
        ,
        \label{Lp_thm_proof:case1_majoration_E(norm(Econd))_partie_x}
    \end{align}
and similarly, 
\begin{align}
    \E \bigcro{\bigp{\int^{y}_{-\infty} \indic_{Y_0 \leq t} \, d\mu(t)}^{1/p}} \leq \int_0^{Q_{(-F_\mu(Y_0))_+^{1/p}}^{-1}((-F_\mu (y))^{1/p})} Q_{(-F_\mu(Y_0))_+^{1/p}} (u) \, du.
    \label{Lp_thm_proof:case1_majoration_E(norm(Econd))_partie_y}
\end{align}

\noindent 
Furthermore, we can select $x>0$ such that $F_\mu(x)^{1/p} = Q_{(F_\mu(Y_0))_+^{1/p} }(\Tilde \beta_{1,Y}(n))$. Indeed, assume that such an $x$ doesn't exist. Then, there exists $x_1>0$ such that $F_\mu(x_1^-) \neq F_\mu(x_1)$ and $ Q_{(F_\mu(Y_0))_+^{1/p} }(\Tilde \beta_{1,Y}(n))$ belongs to $] F_\mu(x_1^-)^{1/p} , F_\mu(x_1)^{1/p} [$. Note that $(F_\mu(Y_0))_+^{1/p}$ doesn't take values in $] F_\mu(x_1^-)^{1/p} , F_\mu(x_1)^{1/p} [$. Since $(F_\mu(Y_0))_+^{1/p}$ and $Q_{(F_\mu(Y_0))_+^{1/p}}(U)$ have the same distribution, it implies that $Q_{(F_\mu(Y_0))_+^{1/p}}$ doesn't take values in $] F_\mu(x_1^-)^{1/p} , F_\mu(x_1)^{1/p} [$. There is a contradiction. 
\\
In the same manner, reasoning on $(-F_\mu(Y_0))_+^{1/p}$ , we can select $y <0$ such that $(-F_\mu(y))^{1/p} = Q_{(-F_\mu(Y_0))_+^{1/p} }(\Tilde \beta_{1,Y}(n))$. 
Hence, selecting $x >0$ such that $F_\mu(x)^{1/p} = Q_{(F_\mu(Y_0))_+^{1/p} }(\Tilde \beta_{1,Y}(n))$ and $y <0$ such that $(-F_\mu(y))^{1/p} = Q_{(-F_\mu(Y_0))_+^{1/p} }(\Tilde \beta_{1,Y}(n))$, by combining (\ref{Lp_thm_proof:case1_majoration_E(norm(Econd))_partie_x}) and (\ref{Lp_thm_proof:case1_majoration_E(norm(Econd))_partie_y}) we get
\begin{align}
    &\E \bigcro{ \bigp{\int^{+\infty}_x \indic_{Y_0 > t}\, d\mu(t)}^{1/p}} 
        \leq \int_0^{\Tilde \beta_{1,Y}(n)} Q_{(F_\mu(Y_0))_+^{1/p}}(u) \, du
        \leq \int_0^{\Tilde \beta_{1,Y}(n)} Q_{Y_{p,\mu}}(u) \, du
    \nonumber \\
    \textrm{ and }& \label{Lp_thm_proof:case1_majoration_E(norm(Econd))_parties_x,y}
    \\
    &\E \bigcro{\bigp{\int^{y}_{-\infty} \indic_{Y_0 \leq t} \, d\mu(t)}^{1/p}} 
        \leq \int_0^{\Tilde \beta_{1,Y}(n)} Q_{(-F_\mu(Y_0))_+^{1/p}}(u) \, du
        \leq \int_0^{\Tilde \beta_{1,Y}(n)} Q_{Y_{p,\mu}}(u) \, du
        .
    \nonumber 
\end{align}

\noindent 
On another hand, as $Q_{Y_{p,\mu}}$ is a nonincreasing function, we have 
\begin{align}\label{Lp_thm_proof:majoration_uQ(u)}
   \Tilde \beta_{1,Y}(n) Q_{Y_{p,\mu}}(\Tilde \beta_{1,Y}(n)) 
    \leq \int_0^{\Tilde \beta_{1,Y}(n)} Q_{Y_{p,\mu}} (u) \, du
   .
\end{align}
Starting from (\ref{Lp_thm_proof:case1_majoration_E(norm(Econd))}) and taking into account (\ref{Lp_thm_proof:case1_majoration_E(norm(Econd))_parties_x,y}) and (\ref{Lp_thm_proof:majoration_uQ(u)}), we get 
\begin{align}\label{Lp_thm_proof:case1_majoration_E(norm(Econd))_final}
    \E \norm{\condEp{X_n}{\calF_0}}_{p,\mu} 
    \leq 6 \int_0^{\Tilde \beta_{1,Y}(n)} Q_{Y_{p,\mu}} (u) \, du.
\end{align}
Using the upper bound of (\ref{Lp_thm_proof:case1_majoration_E(norm(Econd))_final}) in (\ref{Lp_thm_proof:ineq_ppal}) with $a = 6$, we derive 
\begin{align}
    \int_0^{\E\norm{\condEp{X_n}{\calF_0}}_{p,\mu}} Q_{\norm{X_0}_{p,\mu}}& \circ G_{\norm{X_0}_{p,\mu}} (u) \, du \nonumber
    \leq 
    24 \int_0^{\Tilde \beta_{1,Y}(n)} Q_{Y_{p,\mu}}^2(u) \, du
    .
    \nonumber
\end{align}
Thus, as soon as (\ref{Lp_thm:H}) is verified, the condition \textit{\ref{cor_cond_suff:item2}} in Corollary \ref{cor_cond_suff} holds and Theorem \ref{thm_ppal} applies. 
  \\
    
    \noindent\textbf{Acknowledgements. }I would like to thank my two advisors J. Dedecker and F. Merlevède for helpful discussions. I also thank the reviewer for his careful reading of the paper and his comments which improved the presentation of the paper.

    \printbibliography

\end{document}